\documentclass[final]{siamltex}
\usepackage{amssymb}
\font\tenBbb=msbm10
\def\Z{\hbox{\tenBbb Z}}
\def\T{{\hbox{\tenBbb T}}}
\font\tenBbb=msbm10

\font\tenBbb=msbm10

\font\tenBbb=msbm10
\newtheorem{cor}{Corollary}[section]
\newtheorem{rem}{Remark}[section]
\newcommand{\R}{{ I\!\!R}}

\def\Dcal{{\cal D}}

\renewcommand{\Re}{{\mbox{ Re}\,}}
\newcommand{\We}{{\mbox{ We}}}

\renewcommand{\o}{{\omega}}
\newcommand{\barr}{\begin{array}}
\newcommand{\earr}{\end{array}}

\renewcommand{\div}{\mbox{div}\,}

\newcommand{\loc}{\mbox{\tiny loc}}
\newcommand{\ga}{\mbox{\bf g}}

\newcommand{\DD}{\mbox{\bf D}}
\newcommand{\WW}{\mbox{\bf W}}

\newcommand{\tendsto}[1]{\renewcommand{\arraystretch}{0.5}
\begin{array}[t]{c}
\longrightarrow \\
{ \scriptstyle #1 }
\end{array}
\renewcommand{\arraystretch}{1}}
\newcommand{\re}[1]{(\ref{#1})}

\title{Newtonian limit for  weakly viscoelastic fluid flows
 of Olroyds' type\thanks{This work was partially
supported by a "contrat CEDRE'' (Universit\'e Paris 13 -
Universit\'e Libanaise).}}

\author{Luc MOLINET \thanks{
 L.A.G.A., Institut Galil\'ee, Universit\'e Paris 13,
 93430 Villetaneuse, France ({\tt molinet@math.univ-paris13.fr})}
        \and Raafat TALHOUK \thanks{Math\'ematiques, Facult\'e des sciences I, Universit\'e Libanaise,
  Beyrouth, Liban ({\tt rtalhouk@ul.edu.lb}). The author was partially
 supported by  "`le programme d'appui des projets de recherche de l'universit\'e Libanaise"
 }}

\begin{document}

\maketitle



\pagestyle{myheadings}
\thispagestyle{plain}
\markboth{L. MOLINET AND R. TALHOUK}{NEWTONIAN LIMIT FOR VISCOELASTIC FLUID FLOWS}

\begin{abstract}
 This paper is concerned with regular flows of incompressible weakly viscoelastic fluids which obey  a
 differential constitutive law of Oldroyd type.
  We study the newtonian limit for weakly viscoelastic fluid
 flows in  $\R^N$ or $\T^N$ for $N=2,\, 3$, when the Weissenberg number (relaxation time measuring the elasticity
 effect in the fluid) tends to zero.
  More precisely, we prove that the velocity field and the extra-stress tensor converge
   in their existence spaces (we examine the Sobolev-$H^s$ theory and the Besov-$B^{s,1}_2$ theory to
 reach the critical case $s= N/2$) to the corresponding newtonian quantities.
  This convergence results are established in the case of ``ill-prepared"' data.
   We deduce, in the two-dimensional case, a new result concerning the global existence
    of weakly  viscoelastic fluids flow. Our approach  makes use of  essentially  two
     ingredients :
    the stability of the null solution  of the viscoelastic fluids flow and the damping effect,
     on the difference between the extra-stress tensor and the tensor of  rate of deformation, induced by the constitutive law of the fluid.

\end{abstract}

\section{Introduction, main results and notations}
In this paper we investigate the  Newtonian limit of weakly
 viscoelastic fluid flows of Oldroyds'type in $  \Omega=\R^N $ or $ \Omega=\T^N $.

 The dynamics of an homogeneous, isothermic and incompressible fluid flows, is described by the  partial differential derivatives system given by :
\begin{equation}\label{pf}
 \left\{
\begin{array}{rcl}
 \rho\Bigl( u'+(u.\nabla)u\Bigr)  & = & f +\div \sigma  \\
\div u & = & 0
\end{array}
\right.
\end{equation}
Here $ \rho>0 $  is the (constant) density and $ f$ is the  external density body forces.
   $u=u(t,x)$ is the velocity vector field and $\sigma=\sigma(t,x)$ is the symmetric stress-tensor,
    which is split into two parts : $\sigma=-pId +\tau$ where $-pId$ is the spherical part ($p=p(t,x)$ the hydrodynamics pressure) and $\tau$ is the tangential part or the extra-stress tensor. The fluid is called Newtonian if $\tau$ can be expressed linearly in terms of the rate of strain tensor $\DD[u]=\frac{1}{2}(\nabla u+\nabla u^{T}),$ i.e.
\begin{equation}\label{l1}
\tau=2\eta \DD[u]
\end{equation}
where $\eta$ is the viscosity coefficient of the fluid (in this case (\ref{pf}) is the Navier-Stokes system). A fluid for which (\ref{l1}) is not valid is called non-Newtonian or complex fluid.

Infortunately no universal constitutive law exists for
non-Newtonian fluids (see for instance \cite{RHN}). In this paper
we consider a class of fluids with memory. For this kind of
fluids, the extra-stress tensor at a time $t$ depends on $\DD[u]$
and its history. A model taking into account this properties is
the Oldroyd's one. The constitutive law of Oldroyd's type
\cite{OL} is given by :
\begin{equation}\label{OL}
\tau +\lambda_{1} \frac{\Dcal_{a}\tau}{\Dcal t} = 2\eta \left(\DD + \lambda_{2}\frac{\Dcal_{a}\DD}{\Dcal t}\right)
\end{equation}
where $0\leq\lambda_{2}<\lambda_{1}$,
$\lambda_{1}$ is the relaxation time and  $\lambda_{2}$ the retardation time. The symbol $\frac{\Dcal_{a}}{\Dcal t}$ denotes an objective (frame indifferent) tensor derivative (see \cite{RHN}). More precisely,
$$
\frac{\Dcal_{a}\tau}{\Dcal t}=\tau' + (u\cdot\nabla)\tau + \tau\WW-\WW\tau -a(\DD\tau+\tau\DD)
$$
with $\WW[u]=\frac{1}{2} (\nabla u -\nabla u^{T})$ is the vorticity tensor and $ a $ is
 a real number verifying $ -1 \le a \le 1 $.
The limit case $\lambda_{1}>0$ and $\lambda_{2}=0$ corresponds to a purely elastic fluid (which is excluded in our analysis), while the limit case $\lambda_{1}=\lambda_{2}=0$ corresponds to a viscous Newtonian fluid.

The constitutive law (\ref{OL}) is not under an evolution form. This equation can be transformed into a  transport equation by spliting the extra-stress tensor into two parts $\tau_{s}+\tau_{p}$, where $\tau_{s}$ corresponds to a Newtonian part (the solvant) and $\tau_{p}$ to the elastic part (the polymer). Setting  $\tau_{s}=2\eta(1-\omega)\DD[u]$, with $\omega$  defined by $0\leq \omega =1-\frac{\lambda_{2}}{\lambda_{1}}\leq 1$, it follows from \re{OL} that
  $\tau_{p}$ satisfies the following transport equation :
\begin{equation}\label{TOL}
\tau_{p} +\lambda_{1} \frac{\Dcal_{a}\tau_{p}}{\Dcal t} = 2\eta \omega\DD[u]\; .
\end{equation}
From now on we shall denote $\tau_{p}=\tau$ and rewrite (\ref{pf}) and (\ref{TOL}) by using dimensionless variables, we obtain the following partial differential system :

\begin{equation}
\label{nonnewtonian}
\left\{
\begin{array}{rcl}
\Re (u'+(u.\nabla)u) -(1-\o)\Delta u +\nabla p & = & f +\div \tau  \\
\div u & = & 0 \\
\varepsilon\Bigl(\tau'+(u.\nabla)\tau + \ga (\nabla u, \tau) \Bigr) +
\tau & = &  2\o \DD [u]
\end{array}
\right. \quad \mbox{ in } \Omega,
\end{equation}
where $ \ga $ is a bilinear tensor valued mapping defined by
$$
\ga (\nabla u, \tau) =\tau \WW[u] - \WW[u] \tau -a(\DD[u]\tau+\tau \;
\DD[u]) ,
$$
$ \Re =\rho \frac{UL}{\eta} $ and $ \varepsilon=\lambda_1 \frac{U}{L} $ are
  respectively the well-known Reynolds number and the Weissenberg number ($ U $
 and $ L $ represent a typical velocity and typical length of the flow).
  It is worth noticing that the Weissenberg number
   is usually  denoted by $ \We $. Here, since we will make
    the Weissenberg number  tend to zero, we prefer to denote it by $ \varepsilon $.
It is crucial to note  that when $ \varepsilon=0 $,
\re{nonnewtonian}
    reduces to the incompressible Navier-Stokes system :
    \begin{equation}
\label{newtonian} \left\{
\begin{array}{rcl}
\Re (v'+(v.\nabla)v) -\Delta v +\nabla p & = & f   \\
\div v & = & 0
\end{array}
\right. \quad \mbox{ in } \Omega \, .
\end{equation}
  On the other hand, from the definition of the retardation parameter
    we observe that $ \omega=1-\mu/\varepsilon $ where $ 0\le \mu< \varepsilon $ is given by
     $ \mu= \lambda_2 \frac{U}{L}$. Therefore the Newtonian limit of \re{nonnewtonian}
      is actually a limit with  two parameters   $ \varepsilon $ and $ \mu $. To simplify
      the study we could drop a parameter by assuming that  the rate $\mu/\varepsilon$
       (or equivalently $ \omega$) is constant as $ \varepsilon $ tends to zero. In this work,
        instead of doing this,
        we will   only
     impose an uniform upper
      bound
 on $\omega\, (=1-\mu/\varepsilon)$ with respect to $ \varepsilon $.
 \\
System \re{nonnewtonian} is completed by the following  initial
conditions
\begin{equation}
 u_{|t=0}=u_0 \mbox{ and }  \tau_{|t=0}=\tau_0 \label{s3}
\end{equation}
Our approach is quite general and uses the two following ingredients :\\
$ \bullet $ The stability of the null solution of
\re{nonnewtonian} for a fixed $ \varepsilon $ (see \cite{CM} on $
\R^N $
 or $ \T^N $ and  \cite{GS}, \cite{FGO},  \cite{LR} for the case of a bounded domain). \\
 $ \bullet $ The damping of factor $ 1/\varepsilon $ on the quantity $ \tau -2\omega \DD[u] $ induced by
  equation \re{nonnewtonian}$_3 $. \\
  Our results in the Sobolev spaces  are valid for $ \Omega= \R^N $ or $ \T^N $ but to simplify
  the expository we will  only consider  $ \Omega=\R^N $ and  give
  the necessary modification to handle the periodic case. \\
  The main idea is to cut $ u $ and $ \tau $ in low and high frequencies at a level depending
  on $ 1/\varepsilon $.
   Roughly speaking, forgotting the nonlinear terms, the high frequency part
   of $ u-v$ ($v$ is the Newtonian solution, see \re{newtonian} associated with the
   initial data $ u_0 $) will satisfy the homogeneous system linearized around the  null solution
    plus a non homogeneous part containing
    a high frequency term of $ v$. But by the Lebesgue monotone convergence
     theorem, this term will tend to zero in the appropriate norms. The stability of the null solution
      (cf. \cite{CM}, \cite{GS}) will then force the high frequencies of $ u-v $ and
       $ \varepsilon^{1/2} \tau $ to remain small (recall that $ (u-v)(0)=0 $). On the other hand,
         the remaining
     frequencies will tend to zero due to the damping effect on $ \tau-2\omega
     \DD[u] $ which we will use in the same time as a smoothing effect.
      We will describe the main steps of the proof  in Section 1.3.\\
 Note that our analysis is in the spirit of  numerous works on the incompressible limit of compressible
 Navier-Stokes equations (see for instance  \cite{Danch} and references therein). However,
  our analysis is in some aspects easier since there is a damping effect relating to the
   small parameter whereas in the incompressible limit it is   a dispersive effect.\\

 In our knowledge, no such result exists in the literature concerning our study, i.e.  the newtonian limit of non-newtonian
  fluid flows. Moreover,
  our global existence result for regular weakly viscoelastic fluids flow in dimension 2 (see
   Corollary \ref{coroHs}) is new and, in particular,
   not contained in the global existence results of \cite{CM}.
 \subsection{Function spaces and notations.}
In the sequel $C$ denotes a positive constant which may differ at each appearance.  When writing $x\backsimeq y$ (for
$x$ and $y$ two non negative real numbers), we mean that there exist $C_1$ and $C_2$ two positive constants (which do
not depend of $x$ and $y$) such that $C_1 x\leq y \leq C_2x$.  When writing $x\lesssim y$ (for $x$ and $y$ two non
negative real numbers), we mean that there exists $C_1$ a positive constant (which does not depend of $x$ and $y$) such
that $ x\leq C_1 y$. \\
${\cal P} $ will denote the Leray projector on solenoidal vector fields. \\
For $ 1\le p,q \le \infty $, we denote by $ \|\cdot \|_{L^p} $ the usual Lebesgue norm on
 $ \Omega=\R^N $,
 $$
\|v\|_{L^p} =\Bigl( \int_{\R^N} |v|^p(x) \, dx\Bigr)^{1/p} \quad
 $$
and  by  $ \|\cdot \|_{L^q_t L^p} $ the space-time Lebesgue norm on $ ]0,t[ \times \Omega $,
$$
\|v\|_{L^q_t L^p} =\Bigl[ \int_0^t  \|v(\tau)\|^q_{L^p}  \, d\tau\Bigr]^{1/q} \quad
$$
with the obvious modification for $ p,q=\infty $.
For $ s\in \R $, we denote by $ \|\cdot \|_{H^s} $ the usual Sobolev norms on  $ \Omega=\R^N $,
$$
\|v\|_{H^s}=\Bigl(  \int_{\R^N} (1+|\xi|^2)^{s} |\hat{v}(\xi)|^2 \, d\xi \Bigr)^{1/2}
$$
where $ \hat{v} $ is the Fourier transform of $ v$.
 The corresponding scalar product will be denoted by  $ ((\cdot, \cdot ))_{H^s} $ .
  Finally, for any $ \varepsilon>0 $ we introduce the following Fourier projectors
 \begin{equation}
 \widehat{P_\varepsilon f}(\xi) = \chi_{[0,\varepsilon^\alpha]}(|\xi|)\hat{f}(\xi)\quad  \mbox{ and }\quad
 \widehat{Q_\varepsilon f}(\xi) = \chi_{]\varepsilon^\alpha,\infty[}(|\xi|)\hat{f}(\xi)\; ,
 \label{projectors}
 \end{equation}
 where $ \alpha>0 $ will be specified later.
\subsubsection{Homogeneous Besov spaces}
Let $\psi$ in ${\cal{S} }(\R)$  such that $\hat{\psi}$ is supported by the set $\{ z \,
/ \, 2^{-1} \leq |z |\leq 2 \,\}$ and such that
\begin{equation}
  \sum_{j\in \Z} \hat{\psi }(2^{-j}z)=1  \; , \; z \neq 0 \;. \label{psi}
\end{equation}
Define $\varphi$ by
\begin{equation}
  \hat{\varphi}=1-    \sum_{j \geq 1} \hat{\psi}(2^{-j}z) \;,
\end{equation}
and observe that $\varphi \in {\cal{D}} (\R)$, $\hat{\varphi }$ is supported by the ball $\{ z \, / \, |z |\leq 2
\,\}$ and $\hat{\varphi}=1$ for $|z |\leq 1$. We denote now by $\Delta_j$ and $S_j$ the convolution operators on $ \R^N $ whose
symbols are respectively given by $\hat{\psi}(2^{-j}|\xi|)$ and $\hat{\varphi}(2^{-j}|\xi|)$
 where $ \xi\in \R^N $ and $|\xi|=\sqrt{\xi_1^2+..+\xi_N^2}$. Also we define the operator
$\tilde{\Delta}_j$ by
$$
 \tilde{\Delta}_j=\Delta_{j-1}+\Delta_j+\Delta_{j+1} \;,
$$
which satisfies,
$$
 \tilde{\Delta}_j \circ \Delta_j=\Delta_j \;.
$$
\vskip0.3cm For $s$ in $\R$, the homogenous Besov space $ {B}^{s,1}_2(\R^N)$
 (to simplify the notation we will simply denoted it by $ B^{s} ( \R^N)$ ) is the completion of
${\cal{S}}(\R^N)$ with respect to the semi-norm
\begin{equation}
  \| f \|_{B^{s}}=\left\| \{ 2^{js}\| \Delta_j (f) \|_{L^2} \} \right\|_{l^1(\Z)} \; .
\end{equation}
\subsection{Main results}
\begin{theorem} \label{theoHs}
Let $ N=2,3 $ and let $ (u_0,\tau_0)\in H^s(\R^N)\times H^s(\R^{N^2}) $ and $ f\in L^2_{\mbox{loc}}(\R;H^{s-1}) $
 with $ s>N/2 $.  Let $v $ be the Newtonian solution satisfying \re{newtonian} with initial data $ u_0 $ and let
  $ 0<T_0\le \infty $ such that $v\in C([0,T_0];H^s) $.
Then for any $ \delta\in ]0,1[ $ there exists
 $$ \varepsilon_0=\varepsilon_0(N,\Re, \delta, \|v\|_{L^\infty_{T_0} H^{s}},
 \|\nabla v\|_{L^2_{T_0} H^{s}},\|\tau_0\|_{H^s},\|{\cal P} f\|_{L^2_{T_0} H^{s-1}})>0
 $$
  such that for any
  $ 0<\varepsilon<\varepsilon_0 $ the system \re{nonnewtonian}, with
  \begin{equation}
0< \omega \le 1-\delta, \label{ioi}
  \end{equation}
  admits a unique solution
  $$
  u_\varepsilon \in C([0,{T_0}];H^s), \quad \nabla u_\varepsilon
   \in  L^2(0,{T_0};H^{s}) , \; \tau_\varepsilon \in C([0,{T_0}]; H^s) \; .
  $$
Moreover,
\begin{equation}
u_\varepsilon \tendsto{\varepsilon\to 0} v \mbox{ in } C([0,{T_0}];
H^s)\label{convuHs}\, ,
\end{equation}
\begin{equation}
\tau_\varepsilon - 2\omega
\DD[u_\varepsilon]  \tendsto{\varepsilon\to 0} 0 \mbox{ in } L^2(0,{T_0};H^s) \label{convtauHs} \, ,
\end{equation}
\begin{equation}
\varepsilon^{1/2} \tau_\varepsilon \tendsto{\varepsilon\to 0} 0
 \mbox{ in } C([0,{T_0}];H^s) \, .\label{epsilontau}
\end{equation}
\end{theorem}

Recalling that in dimension two, the solution of the Newtonian
problem exists for all positive times, we  deduce the following
result.
\begin{cor} \label{coroHs}
In dimension 2  there exists  $$  \varepsilon_0=\varepsilon_0(\Re, \delta,
\|v\|_{L^\infty_\infty H^{s}},
 \|\nabla v\|_{L^2_{\infty} H^{s}},\|\tau_0\|_{H^s},\|{\cal P} f\|_{L^2_{\infty} H^{s-1}})>0 $$
  such that for any
  $ 0<\varepsilon<\varepsilon_0 $ the solution of \re{nonnewtonian} given by Theorem \ref{theoHs}
  exists for all positive times.
\end{cor}

\begin{rem} Note that Theorem \ref{theoHs} is a convergence result
for ``ill-prepared" data. Indeed the quantity $ \tau_0-2\omega
\DD[u_0] $ is not assumed to be small with $ \varepsilon $.
Moreover, this is a singular limit result since $ \tau $ and $
\DD[u] $ do not belong to the same function space. In particular,
$ \DD[u_0] $ is not as the same level of Sobolev regularity as $
\tau_0 $.
\end{rem}

\begin{rem}
According to the introduction, the Newtonian limit process is actually a limit process with two parameters $ \varepsilon $ and $ \mu $ tending to zero with $ 0\le\mu< \varepsilon $.
The assumption \re{ioi} of Theorem \ref{theoHs} means that we impose the following
  additional conditions on the rate $ \mu/\varepsilon $ as $(\varepsilon,\mu) $
   tends to zero $(0,0) $:
 $$
\delta \le\frac{\mu}{\varepsilon}=\frac{\lambda_2}{\lambda_1}< 1
$$
for some  fixed $ 1>\delta>0  $.
\end{rem}

As mentioned in the introduction, the use of Besov spaces enables us to reach the critical
 index $ s=N/2 $.
\begin{theorem} \label{theoBs}
  Let $ N=2,3 $ and let $ (u_0,\tau_0)\in
B^{N/2-1}(\R^N)\times B^{N/2}(\R^{N^2}) $ and $
f\in L^1_{\mbox{loc}}( B^{N/2-1}) $.  Let $v $ be the Newtonian solution
 satisfying \re{newtonian} with initial data $ u_0 $ and let
  $ 0<T_0\le \infty $ such that $v\in C([0,T_0]; B^{N/2-1}) $.
There exist $ 0<\omega_0<1 $ and
$
\varepsilon_0=\varepsilon_0(N,\Re,\omega_0,\|\tau_0\|_{B^{N/2}},{\cal P} f , u_0) >0 $ such that for any
  $ 0<\varepsilon<\varepsilon_0 $ the system \re{nonnewtonian},
  with
  $
0< \omega\le \omega_0
  $,
  admits a unique solution
  $$
  u_\varepsilon \in C([0,{T_0}]; B^{N/2-1}), \quad  u_\varepsilon
   \in  L^1(0,{T_0}; B^{N/2+1}) , \; \tau_\varepsilon \in C([0,{T_0}];  B^{N/2}) \; .
  $$
Moreover,
\begin{equation}
u_\varepsilon \tendsto{\varepsilon\to 0} v \mbox{ in }
C([0,{T_0}];  B^{N/2-1})\label{convuBs}\; ,
\end{equation}
\begin{equation}
\tau_\varepsilon - 2\omega
\DD[u_\varepsilon] \tendsto{\varepsilon\to 0} 0  \mbox{ in } L^1(0,{T_0}; B^{N/2})
\label{convtauBs} \; ,
\end{equation}
\begin{equation}
\varepsilon^{1/2} \tau_\varepsilon \tendsto{\varepsilon\to 0} 0
 \mbox{ in } C([0,{T_0}]; B^{N/2}) \, .\label{epsilontauBs}
\end{equation}
\end{theorem}

In dimension two, using  the classical global existence result in $ B^{0} $ for the
Newtonian problem  (see for instance
\cite{Danch}), we get the following corollary.
\begin{cor} \label{coroBs}
In dimension 2  there exists  $$ \varepsilon_0=\varepsilon_0(\Re,
\omega_0,\|v\|_{L^\infty_\infty B^{0}},
 \|\nabla v\|_{L^1_\infty B^{2}},\|\tau_0\|_{B^{1}},\|{\cal P} f\|_{L^1_{T_0} B^{0}}) $$
  such that for any
  $ 0<\varepsilon<\varepsilon_0 $ the solution of \re{nonnewtonian} given by Theorem \ref{theoBs}
  exists for all positive times.
\end{cor}
In the surcritical case, $s>N/2$,  we get similar results by considering  non homogeneous Besov spaces. Note that  $ \varepsilon_0$  depends then explicitly on  some norms of $ v$.
\begin{theorem} \label{critic}
For $ s>N/2 $, Theorem \ref{theoBs} and Corollary \ref{coroBs}
also hold by replacing the function spaces $ B^{N/2-1} $  by $
B^{s-1}\cap
 B^{N/2-1} $ and  $ B^{N/2} $ by $ B^{s}\cap
 B^{N/2} $. Moreover, $ \varepsilon_0 $ will depend now explicitly on
  some norms of $ v $ and $ {\cal P} f $.
  More precisely, for $ s>N/2 $, we have
  $$  \varepsilon_0=\varepsilon_0(N,\Re, \omega_0,\|v\|_{L^\infty_{T_0} B^{N/2-1}},
 \|\nabla v\|_{L^1_{T_0} B^{s}},\|\tau_0\|_{B^{N/2}},\|{\cal P} f\|_{L^1_{T_0} B^{N/2-1}})
 \, .
 $$
\end{theorem}

\subsection{Sketch of the proof of Theorem \ref{theoHs}}
In this subsection we want to explain the main steps of the proof
of Theorem \ref{theoHs}. Note that Theorems \ref{theoBs}  follows from the same arguments. To simplify we drop the nonlinear terms in \re{nonnewtonian}. The first step consists in noticing that $W:=u-v $ satisfies the following system:
\begin{equation}
\label{nonnewtonianlin} \left\{
\begin{array}{l}
\Re W_t -(1-\o)Q_\varepsilon \Delta W-P_\varepsilon \Delta W =
  P_\varepsilon {\cal P} (\div \tau -\omega\Delta u)  \\
  \hspace*{3cm} -\omega Q_\varepsilon \Delta v
  + Q_\varepsilon {\cal P}\div \tau  \\
\div W  =  0 \\
\varepsilon\tau_t +Q_\varepsilon \tau  =  2\omega Q_\varepsilon \DD[W] +2\omega
Q_\varepsilon \DD [v] -P_\varepsilon (\tau-2\omega \DD[u])
\end{array}
\right.
\end{equation}
where $ P_\varepsilon $ and $ Q_\varepsilon $ are the projectors on respectively the low and the high frequencies defined in \re{projectors}. \\
Projecting on the high frequencies with $ Q_\varepsilon $ (see
 \re{projectors} for the definition), proceeding as
in \cite{CM}, it is easy to check that we get a
differential inequality close to
\begin{eqnarray*}
  \frac{d}{dt} \Bigl( \|Q_{\varepsilon}W\|_{H^s}^2& +& \varepsilon \| Q_{\varepsilon}\tau\|_{H^s}^2\Bigr)
  + \|Q_{ \varepsilon}\nabla W\|_{H^s}^2
   + \|Q_\varepsilon
   \tau\|_{H^s}^2\lesssim  \|Q_\varepsilon \nabla v\|_{H^s}^2
\end{eqnarray*}
where we drop all the constant to clarify the presentation.
Therefore, since $ W(0)=0 $, $ \varepsilon\to 0 $ and, by the
Lebesgue monotone convergence theorem, $  \|Q_\varepsilon \nabla
v\|_{L^2_{T_0} H^s} \to 0 $, we infer that $ \|Q_{\varepsilon}W\|_{L^\infty_{T_0} H^s} $ goes to zero
with $ \varepsilon $. Now, to treat the low frequency part, we
observe that, computing $ P_\varepsilon (\re{nonnewtonianlin}_3 -
\frac{2\omega}{\Re} \DD[\re{nonnewtonianlin}_1] $ and taking the $
H^{s-1}$-scalar product of the resulting equation with
$Z:=\tau-2\omega \DD[u] $, we obtain  something like
\begin{equation}
  \frac{d}{dt}\|P_{\varepsilon}Z\|_{H^{s-1}}^2 +
  \frac{1}{\varepsilon} \| P_{\varepsilon}Z \|_{H^{s-1}}^2
 \lesssim  \|P_\varepsilon \tau \|_{H^s}^2+\|P_\varepsilon f
 \|_{H^{s-1}}^2
 \label{ru} \; .
\end{equation}
On the other hand, $ P_\varepsilon \re{nonnewtonianlin}_3 $ can be
rewritten as
$$
\varepsilon P_\varepsilon \tau_t
 +\varepsilon^\beta P_\varepsilon \tau =
  2\omega \varepsilon^\beta P_\varepsilon \DD[W]+2\omega
\varepsilon^\beta P_\varepsilon \DD[v]
  -(1-\varepsilon^\beta) P_\varepsilon Z \; ,
$$
where $ 0<\beta<1 $ will be specified later. Therefore, taking
the $ H^s $-scalar product of this last equality with $ \tau $
and adding with the scalar product of \re{nonnewtonianlin}$_1 $ with
 $ W $ we get a differential inequality close to
 $$
  \frac{d}{dt} \Bigl( \|P_{\varepsilon}W\|_{H^s}^2 + \varepsilon
  \| P_{\varepsilon}\tau\|_{H^s}^2\Bigr)
  + \|P_{ \varepsilon}\nabla W\|_{H^s}^2
   + \varepsilon^{\beta} \|P_\varepsilon
   \tau\|_{H^s}^2\lesssim \varepsilon^{-\beta}\|P_\varepsilon Z\|_{H^s}^2
   + \varepsilon^{\beta} \|P_\varepsilon \nabla v\|_{H^s}^2\, .
$$
 Adding this last inequality and $ \varepsilon^{2\beta}  \re{ru}$ we finally
 obtain
 \arraycolsep1pt
  \begin{eqnarray*}
  \frac{d}{dt} \Bigl( \|P_{\varepsilon}W\|_{H^s}^2& +& \varepsilon
  \| P_{\varepsilon}\tau\|_{H^s}^2\Bigr)
   +  \|P_{ \varepsilon}\nabla W\|_{H^s}^2
   + \varepsilon^{\beta} \|P_\varepsilon
   \tau\|_{H^s}^2 +\varepsilon^{2\beta-1}\|P_\varepsilon
   Z\|_{H^{s-1}}^2\\
   & \lesssim  & \varepsilon^{2\beta}\|P_\varepsilon f \|_{H^{s-1}}^2
   + \varepsilon^{\beta} \|P_\varepsilon \nabla v\|_{H^s}^2
\end{eqnarray*}
since $\varepsilon^{-\beta}\|P_\varepsilon Z\|_{H^s}^2 \le
 \varepsilon^{-\beta}\varepsilon^{-2\alpha} \|P_\varepsilon Z\|_{H^{s-1}}^2
\le \frac{\varepsilon^{2\beta-1}}{4}
  \|Z\|_{H^{s-1}}^2 $ as soon as $ 1-3\beta -2\alpha>0 $. This last inequality enables us to
   conclude for the low frequency part. Note that we used the damping effect
  also as a smoothing effect.
\section{Proof of Theorem \ref{theoHs}.}
Let us recall the following existence theorem proven by J.-Y.
Chemin and N. Masmoudi \cite{CM}.
\begin{theorem} \label{CMSobolev} Let $ (u_0,\tau_0) \in H^s(\R^N)\times
H^s(\R^{N^2}) $ with $ s>N/2 $. Then there exists a unique
positive maximal time $ T^* $ and a unique solution
$$
(u,\tau) \in C([0,T^*[;H^s)\cap L^2_{loc}(0,T^*;H^{s+1})\times
C([0,T^*[; H^s)
$$
Moreover, if $ T^*<\infty $ then $ \forall\;  N/2<s'\le s $
\begin{equation}
\limsup_{t\nearrow T^*}
\Bigl(\|u(t)\|_{H^{s'}}+\|\tau(t)\|_{H^{s'}} \Bigr)=+\infty \label{explose}
\end{equation}
\end{theorem}

\begin{rem}
Actually in  \cite{CM} the following sharper blow up condition is derived
 $$
 T^*<\infty \Longrightarrow \int_0^{T^*} \|\nabla u(t)\|_{L^\infty}+\|\tau(t)\|_{L^\infty}^2 \, dt = +\infty \;,
 $$
 but for our purpose the  classical blow-up condition \re{explose} will be sufficient.
\end{rem}

 Let us also recall a  commutator estimate and  classical Leibniz rules for
  fractional derivatives.
\begin{lemma} \label{comm}
Let $\Delta$ be the Laplace operator on $\R^N$, $N\geq 1$. Denote
by $J^{s}$ the operator $(1-\Delta)^{s/2}$. \\
$ \bullet $ For every  $s> N/2 $, \begin{equation}
 \|[J^s,f] g\|_{L^2(\R^N)} \lesssim  \|\nabla
f\|_{H^s(\R^N)} \| g\|_{H^{s-1}(\R^N)} \; .\label{commuHs}
\end{equation}
$ \bullet $ For every $ s>0 $, $ 1<q,q'\le \infty $ and  $
1<r,p,p'<\infty $ with $ 1/p+1/q=1/p'+1/q'=1/r $,
\begin{equation}
\|J^s(fg)\|_{L^r(\R^N)} \lesssim \|J^s f \|_{L^p(\R^N)}
\|g\|_{L^q(\R^N)} +\| f \|_{L^{q'}(\R^N)} \|J^s
g\|_{L^{p'}(\R^N)}\; . \label{prod1}
\end{equation}
$ \bullet $ For every $  p,r,t $ such that $ r,p\ge t $, $ r+p\ge 0 $  and $ r+p-t>N/2 $,
\begin{equation}
\|fg\|_{H^t(\R^N)} \lesssim \|f \|_{H^p(\R^N)} \|g\|_{H^r(\R^N)}
\; . \label{prod2}
\end{equation}
\end{lemma}

{\it Proof. } \re{prod1} and \re{prod2} are classical and can be
found in \cite{KPV4} and \cite{Grisvard}. \re{commuHs} is a variant
of Kato-Ponce's commutator estimates. It is proven in \cite{Tom}
in dimension 1 but the proof works also in dimension 2 and 3.  \vspace*{2mm} \\
To treat some nonlinear terms  in dimension 2  we will need
moreover the following Gagliardo-Nirenberg type inequality (see
for instance \cite{Evans}).
 \begin{lemma} \label{Sobolev2}
 Let $N\geq 2$, for $u\in H^1(\R^N)$ the following Sobolev type inequality hold for any $2\leq p <+\infty $ such that $\frac{1}{2} -\frac{1}{N} \leq \frac{1}{p}$ :
 \begin{equation}
 \|u\|_{L^p(\R^N)}\lesssim \, \|u\|_{L^2(\R^N)}^{(\frac{N}{p}-\frac{N}{2} +1)}
  \|\nabla u \|_{L^2(\R^N)}^{(\frac{N}{2}-\frac{N}{p})} \quad .\label{sobo}
 \end{equation}
 \end{lemma}

\subsection{Estimate on $ W=u-v $ and $ \varepsilon^{1/2} \tau$.}
We start by deriving a differential inequality for the $ H^s
$-norms of $ W $ and $ \varepsilon^{1/2} \tau $. The high
frequency part of this inequality is directly inspired
 by  the stability proof of the
 null solution in \cite{CM}. This will enable us to control the very high frequency
 part $ (Q_\varepsilon u, Q_\varepsilon \tau) $ of the solution. The other part $ (P_\varepsilon u, P_\varepsilon \tau) $
  will be treated by using
  the damping effect. \\
  For $ \varepsilon>0 $ fixed, Theorem \ref{CMSobolev} gives the existence and uniqueness of the solution
   $ (u_\varepsilon, \tau_\varepsilon) $ of \re{nonnewtonian} in
   $ C([0,T_\varepsilon^*[; H^s) \cap L^2_{\loc} (0,T_\varepsilon^*; H^{s+1})\times
    C([0,T_\varepsilon^*[; H^s) $ for some $T_\varepsilon^*>0 $. To simplify the notations, we drop
     the index $ \varepsilon $ on $ u$ and $ \tau $ in the sequel.
  Setting
  $$
  Z= \tau -2\omega \DD[u]
  $$
  we have the following estimates :
  \begin{lemma} For $ \varepsilon>0 $ small enough,
  the solution $ (u,\tau) $ of \re{nonnewtonian} satisfies for all $
0<t<T^*_\varepsilon $  and $ 0<\beta< 1$, \arraycolsep1pt
\begin{eqnarray}
  \frac{d}{dt} \Bigl( \frac{\Re}{2} \|W\|_{H^s}^2& +& \frac{\varepsilon}{4\omega}  \| \tau\|_{H^s}^2\Bigr) +
\frac{1}{4}\|P_\varepsilon\nabla W\|_{H^s}^2
  + \frac{(1-\omega)}{2}\|Q_{ \varepsilon}\nabla W\|_{H^s}^2 \nonumber \\
   &  & +\frac{1}{4\omega} \|Q_\varepsilon
   \tau\|_{H^s}^2+\frac{\varepsilon^\beta}{4\omega} \| P_\varepsilon\tau\|_{H^s}^2\nonumber\\
& \le & (1+\frac{4\varepsilon^{-\beta}}{\omega})\|P_\varepsilon
Z\|_{H^s}^2
  +4\omega \|Q_\varepsilon \nabla v\|_{H^s}^2+8 \omega\varepsilon^\beta\|P_\varepsilon \nabla v\|_{H^s}^2
  \nonumber \\
 & &\hspace*{-20mm}+\frac{C \,
\Re}{(1-\omega)^2} (\|\nabla u\|_{H^s}^2+\|\nabla v\|_{H^{s}}^2 ) \|W\|_{H^s}^2+\frac{C}{\omega} \, \varepsilon^{2-\beta}
  \|\nabla u\|_{H^s}^2
\|\tau\|_{H^s}^2 \label{Wtau}
 \end{eqnarray}
 whenever $ 0<\omega < 1 $. Moreover, for $ 0< \omega \le 10^{-2} $, it holds
 \begin{eqnarray}
  \frac{d}{dt} \Bigl( \frac{\Re}{4} \|W\|_{H^s}^2& +& \frac{\varepsilon}{2}  \| \tau\|_{H^s}^2\Bigr) +
\frac{1}{8}\|P_\varepsilon\nabla W)\|_{H^s}^2
  + \frac{(1-\omega)}{4}\|Q_{ \varepsilon}\nabla W\|_{H^s}^2 \nonumber \\
   &  & +\frac{1}{4} \|Q_\varepsilon
   \tau\|_{H^s}^2+\frac{\varepsilon^\beta}{4} \| P_\varepsilon\tau\|_{H^s}^2\nonumber\\
& \le & (1+4\varepsilon^{-\beta})\|P_\varepsilon
Z\|_{H^s}^2
  +8\omega^2 \|Q_\varepsilon \nabla v\|_{H^s}^2+8 \omega^2\varepsilon^\beta\|P_\varepsilon \nabla v\|_{H^s}^2
  \nonumber \\
 & &\hspace*{-20mm}+\frac{C \,
\Re}{(1-\omega)^2} (\|\nabla u\|_{H^s}^2+\|\nabla v\|_{H^{s}}^2 ) \|W\|_{H^s}^2+C \, \varepsilon^{2-\beta}
  \|\nabla u\|_{H^s}^2
\|\tau\|_{H^s}^2 \label{Wtau22} \; .
 \end{eqnarray}
  \end{lemma}

  {\it Proof. }
Notice that $ W$ verifies the equation
\arraycolsep1pt
\begin{eqnarray}
\Re \Bigl( W_t +{\cal P}(u.\nabla)W\Bigr)&  - & \Delta W  = {\cal P}\div \tau -\omega \Delta u -\Re \, {\cal P}(W.\nabla) v \nonumber\\
& = & P_\varepsilon {\cal P}\Bigl( \div\tau -\omega \Delta u\Bigr) -\omega Q_\varepsilon \Delta v \nonumber\\
 & & + Q_\varepsilon {\cal P} \div \tau
-\omega Q_\varepsilon \Delta W -\Re \, {\cal P}(W.\nabla) v \; .
\label{eqW}
\end{eqnarray}
 Therefore, multiplying scalarly \re{eqW} by $  W $ in $
H^s(\R^N) $, using Cauchy-Schwarz, Lemma \ref{comm} and that $ u $ is divergence free, we obtain
 \begin{eqnarray}
 \frac{1}{2} & \Re &  \frac{d}{dt} \| W\|_{H^s}^2 + \| P_\varepsilon\nabla W\|_{H^s}^2
  + (1-\omega)\|Q_{ \varepsilon}\nabla W\|_{H^s}^2
 \nonumber  \\
& \le &  ((Q_\varepsilon \div \tau, W ))_{H^s} +\|
P_\varepsilon Z\|_{H^s} \|\nabla W\|_{H^s}
  +\omega  \|Q_\varepsilon \nabla v\|_{H^s}  \|Q_\varepsilon\nabla W\|_{H^s} \nonumber \\
   & & +C
\Re \Bigl|((J^s(u.\nabla)W, J^s W))_{L^2} \Bigr|
+\Re \Bigl|((J^s(W.\nabla)v, J^s W))_{L^2} \Bigr| \label{trtr}
 \end{eqnarray}
 To estimate the second to the last term of \re{trtr}, we rewrite it with the help of a commutator and
  apply Cauchy-Schwarz to the term containing this commutator to get
  \begin{eqnarray*}
  \Bigl|((J^s(u.\nabla)W, J^s W))_{L^2} \Bigr| \le   \Bigl|((\,(u.\nabla)J^s W, J^s W))\Bigr|+
\Bigl\|[J^s, (u.\nabla)] W\Bigr\|_{L^2}
  \|J^s W\|_{L^2} \;.
  \end{eqnarray*}
  Since $ u $ is divergence free, the first term of the right-hand side
of this last inequality cancels by integration by parts. Estimating the second term
thanks to Lemma \ref{comm}, we then obtain
 $$
  \Bigl|((J^s(u.\nabla)W, J^s W))_{L^2} \Bigr| \le C \, \|\nabla u\|_{H^s} \|\nabla W\|_{H^s} \|W\|_{H^s}\quad .
 $$
 Now, to estimate the last term of the right-hand side \re{trtr} we have to distinguish the cases $ N=2 $ and
$ N=3 $. \\
$ \bullet $ $ N=3 $. Then by Lemma \ref{comm}, H\"older, Sobolev and Young inequalities, we get
 \begin{eqnarray*}
 \Bigl|((J^s(W.\nabla)v, J^s W))_{L^2} \Bigr| & \le & \|J^s(W.\nabla)v\|_{L^{6/5}} \|J^s W\|_{L^6} \\
 & \lesssim &   \Bigl( \|J^s W\|_{L^2} \|\nabla v \|_{L^3}+\|W\|_{L^3} \|J^s \nabla v\|_{L^2} \Bigr) \|J^s \nabla W\|_{L^2} \\
 & \lesssim &   \|W\|_{H^s} \|\nabla v\|_{H^s} \|\nabla W\|_{H^s}
 \quad .
 \end{eqnarray*}
$ \bullet $ $ N=2 $. In this case, using H\"older and  Lemmas \ref{comm}-\ref{Sobolev2}  we infer that
\begin{eqnarray*}
 \Bigl|((J^s(W.\nabla)v), J^s W))_{L^2} \Bigr| & \le & \|J^s(W.\nabla)v\|_{L^{3/2}} \|J^s W\|_{L^3} \\
 & \lesssim &   \Bigl( \|J^s W\|_{L^6} \|\nabla v \|_{L^2}+\|W\|_{L^6} \|J^s \nabla v\|_{L^2} \Bigr) \|J^s  W\|_{L^3} \\
 & \lesssim &   \|J^s W\|_{L^6} \|J^s \nabla v \|_{L^2}  \|J^s  W\|_{L^3} \\
 & \lesssim &  \|W\|_{H^s} \|\nabla v\|_{H^s} \|\nabla W\|_{H^s}
 \quad .
 \end{eqnarray*}
  By Young inequalities it thus follows from \re{trtr} that
\begin{eqnarray}
 \frac{\Re}{2} \frac{d}{dt} \| W\|_{H^s}^2& +& \frac{3}{4}\| P_\varepsilon\nabla W\|_{H^s}^2
  + \frac{(1-\omega)}{2} \|Q_{> \varepsilon}\nabla W\|_{H^s}^2
  \nonumber \\
&  & \hspace{-25mm}\le   ((Q_\varepsilon \div \tau, W ))_{H^s}
+
\| P_\varepsilon Z\|_{H^s}^2 +\frac{\omega^2}{4(1-\omega)} \|Q_\varepsilon \nabla v\|_{H^s}^2\nonumber \\
  & & \hspace{-25mm}+\frac{C\, \Re}{(1-\omega)^2} \,
 \Bigl( \|\nabla u\|_{H^s}^2+\|\nabla v\|_{H^s}^2\Bigr) \|W\|_{H^s}^2
 \; .  \label{W1}
 \end{eqnarray}
 On the other hand, for $ 0<\beta<1 $, observing that
\begin{eqnarray*}
\tau-2\omega \DD[u] & = & Q_\varepsilon \tau-2\omega Q_\varepsilon( \DD[W]+\DD[v]) \\
&  & +(1-\varepsilon^\beta)P_\varepsilon Z+\varepsilon^\beta \Bigl( P_\varepsilon \tau
  - 2\omega  P_\varepsilon( \DD[W]+\DD[v]) \Bigr)\; ,
\end{eqnarray*}
we deduce from \re{nonnewtonian} that $ \tau $ satisfies the equation
 \begin{eqnarray*}
 \varepsilon \Bigl(\tau_t +(u.\nabla)\tau & + & g(\nabla u,\tau)\Bigr)+Q_\varepsilon \tau
 +\varepsilon^\beta P_\varepsilon \tau = 2\omega  Q_\varepsilon \DD[W]+2\omega  Q_\varepsilon \DD[v]\nonumber\\
  &  &
 +2\omega \varepsilon^\beta P_\varepsilon \DD[W]+2\omega \varepsilon^\beta P_\varepsilon \DD[v]
  -(1-\varepsilon^\beta) P_\varepsilon Z \; .
 \end{eqnarray*}
Taking the $ H^s $ scalar product of this equation with $ \tau$, using Lemma \ref{comm}, Cauchy-Schwarz and
 Young inequalities we get
\begin{eqnarray}
\frac{\varepsilon}{2} \frac{d}{dt} \|\tau\|_{H^s}^2 & + & \frac{1}{2}\|Q_{ \varepsilon}\tau\|_{H^s}^2
+\frac{\varepsilon^\beta}{2} \| P_\varepsilon\tau\|_{H^s}^2\le  2\omega ((Q_\varepsilon \DD[W], \tau))_{H^s} +8\varepsilon^{-\beta} \|P_\varepsilon Z\|_{H^s}^2
  \nonumber \\
 & &+8\omega^2 \|Q_{ \varepsilon} \nabla v\|_{H^s}^2+8 \omega^2\varepsilon^\beta (\|P_\varepsilon \nabla W\|_{H^s}^2+\|P_\varepsilon \nabla v\|_{H^s}^2) \nonumber \\
 & & +C \, \varepsilon^{2-\beta} \|\nabla u\|_{H^s}^2  \|\tau\|_{H^s}^2 \; . \label{tau1}
\end{eqnarray}
We now separate the two cases : \\
$ \bullet $ $\omega \neq 0 $.
 Then, adding \re{W1} and $\re{tau1}/2\omega $ we notice that the first term in
  the right-hand side of \re{W1} and \re{tau1}
 cancel each other and \re{Wtau} follows. This gives \re{Wtau} for $ \varepsilon $ small enough
  since $ \beta>0 $. \\
 $ \bullet $ $ 0< \omega \le 10^{-2} $. Then adding \re{W1}/2+\re{tau1}, estimating the two remaining $ H^s$-scalar products by integration by parts, Cauchy-Schwarz inequality and Young inequality, one obtains \re{Wtau22}

\subsection{Estimate on $Z=\tau-2\omega \DD[u]$}
We will now take advantage of the damping effect on $Z=\tau-2\omega
\DD[u]$.
\begin{lemma}
 The solution $ (u,\tau) $ of \re{nonnewtonian} satisfies for all $ \varepsilon $ small enough
  and $
0<t<T_\varepsilon^* $,
\begin{eqnarray}
\frac{1}{2}\frac{d}{dt} \|Z\|_{H^{s-1}}^2 &  + &
\frac{1}{2\varepsilon}
  \|Z\|_{H^{s-1}}^2 \le
  \frac{4\omega}{\Re(1-\omega)}
 \|{\cal P}f \|_{H^{s-1}}^2+\frac{(1+\omega)^2}{\Re (1-\omega)} \|\tau \|_{H^{s}}^2
  \nonumber \\
   & &+ \frac{4 }{1-\omega}
 (\Re \|\nabla u\|_{H^{s}}^2+\|\tau \|_{H^s}^2)\|u\|_{H^{s}}^2
  \; . \label{Z2}
 \end{eqnarray}
  \end{lemma}

  {\it Proof. }
 We apply $ \frac{2\omega}{\Re} \DD[\cdot] $  to $\re{nonnewtonian}_1 $
 and substract the resulting equation from  $ \re{nonnewtonian}_2 $ to obtain
 \begin{equation}
 Z_t -\frac{(1-\omega)}{\Re} \Delta Z +\frac{1}{\varepsilon} Z = -f_1-f_2  \label{Z1}
 \end{equation}
 where
 $$
 f_1=\frac{2\omega}{\Re} \DD[{\cal P} \div \tau] -\frac{(1-\omega)}{\Re} \Delta \tau
  +\frac{2\omega}{\Re} \DD[{\cal P} f]
 -2\omega \DD[{\cal P}(u.\nabla)u]
 $$
 and
 $$ f_2= {\cal P}(u.\nabla)\tau +g(\nabla u,\tau) \quad .
 $$
 Taking the $ H^{s-1} $-scalar product of \re{Z1} with $ Z $ we get
 \begin{eqnarray}
 \frac{1}{2}\frac{d}{dt} \|Z\|_{H^{s-1}}^2 &  + &
  \frac{(1-\omega)}{4\Re} \|\nabla Z\|_{H^{s-1}}^2 +\frac{1}{\varepsilon}
  \|Z\|_{H^{s-1}}^2 \nonumber \\
  & & \hspace*{-20mm} \le  C \, \frac{(1+\omega)^2}{\Re (1-\omega)} \|\tau \|_{H^{s}}^2
  +\frac{4\omega}{\Re(1-\omega)} \|{\cal P}f \|_{H^{s-1}}^2 \nonumber \\
   & &\hspace*{-20mm} +
 \frac{4 \omega\Re}{1-\omega}
 \|(u.\nabla) u\|_{H^{s-1}}^2 +4 \|(u.\nabla)\tau\|_{H^{s-1}}^2
 + \|\ga (\nabla u,\tau) \|_{H^{s-1}}^2\, , \label{Z3}
 \end{eqnarray}
 where we used that
 \begin{eqnarray*}
\frac{2\omega}{\Re} \Bigl|((\DD[{\cal P} \div \tau], Z))_{H^{s-1}}
\Bigr| & \le&  \, C \frac{2\omega}{\Re} \|\div \tau \|_{H^{s-1}}
\|\nabla Z \|_{H^{s-1}} \\
& \le & \frac{1-\omega}{8\Re} \|\nabla Z \|_{H^{s-1}}^2+C\,
\frac{\omega^2}{\Re(1-\omega)} \|\tau\|_{H^s}^2 \, .
 \end{eqnarray*}
 Finally to control the nonlinear terms we notice that thanks to   \re{prod2}
 ,
    $$
   \|a .\nabla b \|_{H^{s-1}} \lesssim \|a\|_{H^s} \|\nabla b\|_{H^{s-1}}
    $$
    which concludes the proof of \re{Z2}.
 \subsection{Convergence to the Newtonian flow}
 We give here the proof in the case $ 10^{-2} \le \omega \le 1-\delta$. The case $ 0< \omega \le 10^{-2} $
  is simpler and can be handled in the same way by using \re{Wtau22} instead of \re{Wtau}. \\
 Adding \re{Wtau} and $ \varepsilon^{2\beta} \re{Z2} $, we
  obtain for $ \varepsilon $ small enough
\begin{eqnarray}
  \hspace*{-5mm}\frac{d}{dt} \Bigl( \frac{\Re}{2} \|W\|_{H^s}^2 +\frac{\varepsilon}{4\omega}
   & \| &
\tau\|_{H^s}^2+
 \frac{\varepsilon^{2\beta}}{2} \|Z\|_{H^{s-1}}^2 \Bigr)  \nonumber \\
 + \frac{(1-\omega)}{4} \|\nabla W\|_{H^s}^2& +&
\frac{1}{8\omega} (\|Q_\varepsilon
   \tau\|_{H^s}^2+\varepsilon^\beta \|P_\varepsilon
   \tau\|_{H^s}^2) +  \frac{\varepsilon^{2\beta-1}}{4}
  \|Z\|_{H^{s-1}}^2 \nonumber\\
& & \hspace*{-28mm} \le  \quad
  8\omega^2 \|Q_\varepsilon \nabla v\|_{H^s}^2+8\omega^2 \varepsilon^\beta\|P_\varepsilon \nabla v\|_{H^s}^2
  +C \, \varepsilon^{1-\beta}  \|\nabla u\|_{H^s}^2\, \frac{\varepsilon}{4\omega} \|\tau\|_{H^s}^2
  \nonumber \\
 & &\hspace*{-28mm}+\frac{C \,
\Re}{(1-\omega)^2}( \|\nabla v\|_{H^{s}}^2+ \|\nabla u\|_{H^s}^2)
\|W\|_{H^s}^2
\nonumber \\
& &\hspace*{-28mm}  + C \frac{ \varepsilon^{2\beta}}{ (1-\omega)}
\Bigl[ \frac{1}{\Re} \|{\cal P} f \|_{H^{s-1}}^2  +
 (\Re \|\nabla u\|_{H^{s}}^2+\|\tau\|_{H^{s}}^2) \|u \|_{H^s}^2 \Bigl]
\label{Wtau2}
\end{eqnarray}
Here, we used that for $ \varepsilon $ small enough, $$
\frac{(1+\omega)^2}{\Re (1-\omega)}\varepsilon^{2\beta}
\|\tau\|_{H^s}^2 \le \frac{\varepsilon^\beta}{8\omega}
\|\tau\|_{H^s}^2
$$
and
$$
(1+\frac{4\varepsilon^{-\beta}}{\omega})\|P_\varepsilon Z\|_{H^s}^2 \le
 (1+20^2 \varepsilon^{-\beta})\varepsilon^{-2\alpha} \|P_\varepsilon Z\|_{H^{s-1}}^2
\le \frac{\varepsilon^{2\beta-1}}{4}
  \|Z\|_{H^{s-1}}^2 \quad
 $$
 as soon as $ \beta>0 $ and $ 1-3\beta-2\alpha>0 $. \\
  From now on we thus take to simplify $ (\alpha,\beta)=(1/8,1/8) $.  Setting
\begin{eqnarray*}
X_s(t)& =& \frac{\Re}{2} \|W(t)\|_{H^s}^2 + \frac{\varepsilon}{4\omega}
\|\tau(t)\|_{H^s}^2+
 \frac{\varepsilon^{2\beta}}{2} \|Z\|_{H^{s-1}}^2\\
& & \hspace*{-8mm}+\int_0^t \frac{(1-\omega)}{4} \|\nabla
W\|_{H^s}^2 + \frac{1}{8} \|Q_\varepsilon
   \tau\|_{H^s}^2 +  \frac{\varepsilon^\beta}{8} \|P_\varepsilon
   \tau\|_{H^s}^2 +  \frac{\varepsilon^{2\beta-1}}{4}
  \|Z\|_{H^{s-1}}^2 \, ds
\end{eqnarray*}
we  infer that $ X_s$  satisfies the following differential inequality
\begin{eqnarray}
\frac{d}{dt} X_s & \le &8 \omega^2 \|Q_\varepsilon \nabla
v\|_{H^s}^2+8\omega^2 \varepsilon^\beta\|P_\varepsilon \nabla
v\|_{H^s}^2
+ C \frac{\varepsilon^{2\beta}}{\Re (1-\omega)} \|{\cal P}f \|_{H^{s-1}}^2 \nonumber \\
& & \hspace*{-15mm}+C(\Re,\delta)
 \Bigl[\varepsilon^{2\beta} \|\tau\|_{H^s}^2+ \|\nabla u\|_{H^{s}}^2+\|\nabla v\|_{H^{s}}^2
  \Bigr] X_s \nonumber \\
& & \hspace*{-15mm}+C \frac{\varepsilon^\beta}{(1-\omega)}
(\Re\varepsilon^\beta \|\nabla
u\|_{H^{s}}^2+\varepsilon^\beta\|\tau\|_{H^{s}}^2)
 \|v \|_{H^s}^2
\end{eqnarray}
where we rewrite $ u $ as $ W+v $ and use the triangle inequality
when necessary. Hence,  Gronwall inequality leads to
\begin{eqnarray}
X_s(t) & \le & \exp\Bigl[ C(\Re,\delta) \Bigl(\varepsilon^{2\beta}
\|\tau\|_{L^2_t H^s}^2 + \|\nabla u\|_{L^2_t
H^{s}}^2+\|\nabla v\|_{L^2_t H^{s}}^2\Bigr)\Bigr] \nonumber  \\
& & \Bigl[ X_s(0)+8\omega^2 \|Q_\varepsilon \nabla v\|_{L^2_t
H^s}^2 +8\varepsilon^\beta\| \nabla v\|_{L^2_t H^s}^2+ C\frac{\varepsilon^{2\beta}}{\Re \delta} \|{\cal P}f \|_{L^2_t H^{s-1}}^2\nonumber \\
&  &+\frac{\varepsilon^\beta}{\delta}  (\Re
\varepsilon^\beta\|\nabla u\|_{L^2_t H^{s}}^2
+\varepsilon^\beta\|\tau\|_{L^2_t H^{s}}^2)
 \|v \|_{L^\infty_t H^s}^2
\Bigr]\; .
\end{eqnarray}
Rewriting $ u$ as $ v+W $, we finally obtain
\begin{eqnarray}
X_s(t) & \le & \exp\Bigl[ C(\Re,\delta) \Bigl(\|\nabla v\|_{L^2_t
H^{s}}^2+X_s(t)\Bigr)
\Bigr] \nonumber  \\
& & \hspace*{-10mm}\Bigl[ X_s(0)+8 \omega^2 \|Q_\varepsilon \nabla
v\|_{L^2_t H^s}^2+8\varepsilon^\beta\| \nabla v\|_{L^2_t H^s}^2
+ C\frac{\varepsilon^{2\beta}}{\Re\, \delta} \|{\cal P}f \|_{L^2_t H^{s-1}}^2 \nonumber \\
 & &\hspace*{-10mm} +\,  C\frac{(1+\Re) \varepsilon^\beta}{\delta} \|v\|_{L^\infty_t H^s}^2 X_s(t)
 +C \frac{\Re \varepsilon^{2\beta}}{\delta} \|v\|_{L^\infty_t H^s}^2
\|\nabla v\|_{L^2_t H^s}^2
\Bigr] \label{estfinale}
\end{eqnarray}
where
\begin{equation}
X_s(0)=\frac{\varepsilon}{4\omega}  \|\tau_0\|_{H^s}^2+
 \frac{\varepsilon^{2\beta}}{2}
 \|\tau_0-2\omega\DD|u_0]\|_{H^{s-2}}^2 \quad . \label{X0}
 \end{equation}
 Let us now assume that $ T^*_\varepsilon\le  T_0$.
 Since \re{estfinale} holds for any $ N/2<s'<s $, noticing
 that
 $$
\|Q_\varepsilon \nabla v\|_{L^2_t H^{s'}} \le
\varepsilon^{\alpha(s-s')}  \|\nabla v\|_{L^2_t H^s}\quad
,
 $$
we deduce from \re{estfinale}, \re{X0} and the continuity of
 $ t\mapsto X_{s'}(t) $ that
there exists $ \varepsilon_0(s, \|\tau_0\|_{H^s},
\|u_0\|_{H^{s}}, \|\nabla v\|_{L^2_{T_0} H^s},
\|v\|_{L^\infty_{T_0} H^s})>0 $ such that for any $ 0<\varepsilon
<\varepsilon_0 $ and any $ 0<t<T^*_\varepsilon $
\begin{equation}
X_{s'}(t) \le C \,
\varepsilon^{\min(\beta,2\alpha(s-s'))}\le C \varepsilon^{ \min(1,s-s')/8}
\end{equation}
 which contradicts \re{explose} of Theorem \ref{CMSobolev}. This ensures that $
 T^*_\varepsilon >T_0 $. Now, since by  Lebesgue monotone convergence
 theorem
 \begin{equation}
\|Q_\varepsilon \nabla v\|_{L^2_{T_0} H^{s}} \tendsto{\varepsilon
\to 0} 0 \label{tyty}  \quad ,
 \end{equation}
 it follows from \re{estfinale}-\re{X0} that $ X_s(T_0)\to 0 $ as $
 \varepsilon \to 0 $. This proves \re{convuHs} and \re{epsilontau}.
  To prove \re{convtauHs} we observe that from this last limit and \re{tyty},
 $ \|Q_\varepsilon (\tau-2\omega \DD[u]) \|_{L^2_{T_0} H^s} \to 0
 $ and $ \varepsilon^{2\beta-1} \|P_\varepsilon (\tau-2\omega
 \DD[u]) \|^2_{L^2_{T_0} H^{s-1}} \to 0 $. This yields the result by
 Bernstein inequality  since $ 2\beta-1+2\alpha < 0 $.
\subsection{The periodic setting}
Let us give here the modifications needed to handle with the case
$ \Omega= \T^N $, $ N=2,3 $. It is worth noticing that Lemma
\ref{comm} holds also with $ \Omega=\T^N $. On the other hand,
the Sobolev inequality \re{Sobolev2} does not hold for general
functions in $ \T^N $ but holds, for instance, for zero-mean
value functions. Note that if $ f(t) $ has mean value zero for
all time $ t\ge 0 $ then using the invariance by Galilean
transformations, $ u\mapsto u(t,x-z\,t)+z $ with $ z\in \R^3 $,
we can assume that $ u $ has zero mean-value for all time and we
are done. Otherwise, we have only to care about the treatment of
the nonlinear term $ (W.\nabla)v $ in \re{eqW}.
  Denoting by $ \underline{W} $ the $ L^2 $-projection of $ W $ on zero mean-value functions,
   we rewrite $ (W.\nabla)v $ as
  \begin{equation}
(W.\nabla)v=(\underline{W}.\nabla)v + (\int_{\Omega}W) \nabla v
\quad . \label{eee}
\end{equation}
We take the $ H^s $-scalar product of \re{eqW} with $
\underline{W} $ and add with the $ L^2 $-scalar product of
\re{eqW} with $ W $. The $ H^s $-scalar product coming from the
first term of the right-hand side of \re{eee} can be treated as in
$ \R^N $. For the second term, we observe that
 \begin{eqnarray*} \Bigl|\Bigl(\!\Bigl(
(\int_{\Omega}W)\nabla v,\underline{W} \, \Bigr)\!\Bigr)_{H^s}\Bigr| & = &
|\int_{\Omega}W|\, \Bigl|(( \nabla v,\underline{W} \,
))_{H^s}\Bigr| \\
& \lesssim & \|W\|_{L^2} \|\nabla v \|_{H^s} \|\nabla
\underline{W}\|_{H^{s-1}}
\end{eqnarray*}
On the other hand, concerning the $ L^2 $-scalar product we notice
that
\begin{eqnarray*}
\Bigl|\Bigl(\!\Bigl( (W.\nabla)v,W\Bigr)\!\Bigr)_{L^2}\Bigr| & = & \Bigl|\,
\Bigl(\!\Bigl((\underline{W}.\nabla)v,W\, \Bigr)\!\Bigr)_{L^2}+(\int_{\Omega}W)\, \Bigl(\!\Bigl(
\,\nabla v,\underline{W} \, \Bigr)\!\Bigr)_{L^2}\Bigr| \\
 & \lesssim & \|\underline{W}\|_{L^6} \|\nabla v
 \|_{L^2}\|W\|_{L^3}+|\int_{\Omega} W| \|\nabla v \|_{L^2}
 \|\underline{W}\|_{L^2} \\
 & \lesssim & \|W\|_{H^s}\|\nabla W\|_{L^2} \|\nabla v \|_{L^2}
\end{eqnarray*}
  We thus obtain exactly as in \re{W1},
\begin{eqnarray*}
 \frac{\Re}{2} \frac{d}{dt} \Bigl( \| \underline{W}\|_{H^s}^2+\|W\|_{L^2}^2\Bigr)& +& \frac{3}{4}
 \| P_\varepsilon\nabla W\|_{H^s}^2
  + \frac{(1-\omega)}{2} \|Q_{> \varepsilon}\nabla W\|_{H^s}^2
  \nonumber \\
&  & \hspace{-25mm}\le
((Q_\varepsilon\div\tau,\underline{W}))_{H^s} +
\| P_\varepsilon Z\|_{H^s}^2 +\frac{\omega^2}{4(1-\omega)} \|Q_\varepsilon \nabla v\|_{H^s}^2\nonumber \\
  & & \hspace{-25mm}+\frac{C\, \Re}{(1-\omega)^2} \,
 \Bigl( \|\nabla u\|_{H^s}^2+\|\nabla v\|_{H^s}^2\Bigr) \|W\|_{H^s}^2
 \; .  \label{WW1}
 \end{eqnarray*}
The remainder of the analysis   is now exactly the same as in $
\R^N $.
 \section{Proof of Theorems \ref{theoBs} and \ref{critic}} In this section we prove a
convergence result in the Besov spaces $ B^{s-1,1}_2 $, $ s\ge
N/2$. It will require a smallness
 assumption on the retardation parameter $ \omega $ but on the other hand will  enable us to reach the critical
  regularity space for \re{nonnewtonian}. Note that our smallness assumption
   on the retardation parameter is the same as the one in \cite{CM} to get the
    stability of the null solution in such function spaces.\\
Let us recall the following well-posedness result derived in
\cite{CM}.
\begin{theorem} \label{CMBesov} Let $ (u_0,\tau_0) \in B^{s-1}(\R^N)\cap B^{N/2-1}(\R^N)\times
B^s(\R^{N^2})\cap B^{N/2}(\R^{N^2}) $ with $ s\ge N/2 $. Then
there exists a unique positive maximal time $ T^* $ and a unique
solution
$$
(u,\tau) \in C([0,T^*[;B^{s-1}\cap B^{N/2-1})\cap
L^1_{loc}(0,T^*;B^{s+1}\cap B^{N/2+1})\times C([0,T^*[; B^s\cap
B^{N/2})
$$
Moreover, if $ T^*<\infty $ then
\begin{equation}
\limsup_{t\nearrow T^*}
\Bigl(\|u(t)\|_{B^{N/2-1}}+\|\tau(t)\|_{B^{N/2}} \Bigr)=+\infty
\label{exploseBs}
\end{equation}
\end{theorem}

    We will make use of the following classical commutator and product estimates
     (see for instance \cite{CM}, \cite{Danch} and \cite{Taylor})
\begin{lemma} \label{commubesov}
For all $ s\in ]1-N/2, 1+N/2[ $ we have
\begin{equation}
\|\tilde{\Delta}_j [(a.\nabla),\Delta_j]b \|_2 \lesssim
2^{-j(s-1)} \gamma_j \|\nabla a\|_{B^{N/2+1}} \| b\|_{B^{s-1}},\label{commu}
\end{equation}
with $ \|\gamma_j\|_{L^1(\Z)} \lesssim 1 $.\\
For all $s_1,s_2\le N/2 $ with  $ s_1+s_2> 0 $ it holds
\begin{equation}
\| a b \|_{B^{s_1+s_2-N/2}} \lesssim \|a\|_{B^{s_1}}
\|b\|_{B^{s_2}}\, . \label{product}
\end{equation}
\end{lemma}

For any $ \varepsilon>0 $ we divide $ \Z $ into the three
following subsets
$$
I:=\Z_-^*=\{j\in \Z,0< 2^j <1 \} , \quad J_\varepsilon:=\{j\in \Z,
1\le 2^j \le \varepsilon^{-\alpha}\} \mbox{ and }
K_\varepsilon:=\{j\in \Z,  2^j
> \varepsilon^{-\alpha}\}\;
$$
 and for any subset $ N\subset \Z $ we denote by $ \|\cdot \|_{B^s_N} $
the semi-norm
$$
\|u\|_{B^s_N} =\sum_{j\in N} 2^{js} \|\Delta_j u \|_{L^2} \, .
$$
\subsection{Estimate on  $ W $ and $  \varepsilon \tau $}
\begin{lemma}The solution
 $ (u ,\varepsilon \tau) $ of \re{nonnewtonian} satisfies for all $
0<t<T^* $ \arraycolsep1pt
\begin{eqnarray}
\frac{d}{dt} \Bigl( \Re \| W\|_{B^{s-1}} & + &
 4 \varepsilon \|\tau\|_{B^{s}}\Bigr) \nonumber \\
& & + [(1-\omega)/2-16\omega ] \| W\|_{ B^{s+1}_{K_\varepsilon}}
+\| W\|_{ B^{s+1}_{I\cup J_\varepsilon}}
 +2  \| \tau\|_{ B^{s}_{K_\varepsilon}}+2 \varepsilon^\beta
 \|\tau\|_{B^s_{I\cup
J_\varepsilon}} \nonumber \\
& \le & 5 \| Z \|_{ B^{s}_{I\cup J_\varepsilon}}
 + 16\omega  \|v\|_{ B^{s+1}_{K_\varepsilon}}+16\omega  \varepsilon^\beta
  (\|W\|_{B^{s+1}_{I\cup
J_\varepsilon}}
 +\|v\|_{B^{s+1}_{I\cup
J_\varepsilon}})
   \nonumber \\
 & & \hspace*{-20mm} + C \, \varepsilon \mu_1 \|u\|_{B^{N/2+1}} \|\tau\|_{B^s}+ C \,
 (\| u\|_{ B^{N/2+1}}+\| v\|_{ B^{N/2+1}}) \|W\|_{ B^{s-1}}\, .
 \label{utaubesov}
\end{eqnarray}
\end{lemma}

{\it Proof. } Applying $ \Delta_j $ to \re{eqW} we have for $j \in J_\varepsilon $,
\begin{eqnarray}
\Re \Bigl( \partial_t \Delta_j W &  + &  {\cal P} (u.\nabla)\Delta_j
W\Bigr) - (1-\omega) \Delta_j \Delta W  =   -\omega \Delta_j
\Delta v+ \Delta_j  {\cal P} \div \tau
\nonumber \\
&  & +\Re \tilde{\Delta}_j {\cal P}[(u.\nabla),\Delta_j] W +\Re \Delta_j
{\cal P}(W.\nabla)v \label{u1}
\end{eqnarray}
and for $  j\in I $,
\begin{eqnarray}
\Re \Bigl( \partial_t \Delta_j W &  + &  {\cal P} (u.\nabla)\Delta_j
W\Bigr) -  \Delta_j \Delta W  =  \Delta_j Z
\nonumber \\
&  & +\Re \tilde{\Delta}_j {\cal P}[(u.\nabla),\Delta_j] W +\Re \Delta_j
{\cal P}(W.\nabla)v \, . \label{u2}
\end{eqnarray}
Taking the scalar
product in $ L^2(\R^N) $ of \re{u1} with $\Delta_j W $, using
that $ W$ is divergence free and   Cauchy-Schwarz inequality a we get
\begin{eqnarray}
\frac{1}{2}\Re \frac{d}{dt} \|\Delta_j W\|_2^2 & + & (1-\omega) \|
\nabla \Delta_j  W\|_2^2 \le \|\Delta_j W\|_2 \Bigl( \omega \|\Delta_j\Delta v\|_2+
\|\Delta_j \div \tau\|_2 \nonumber\\
 & & \hspace*{-4mm} +\Re \|\tilde{\Delta}_j P[(u.\nabla),\Delta_j] W\|_2 +\Re \|\Delta_j
P(W.\nabla)v\|_2 \Bigr). \label{n2}
\end{eqnarray}
We use now that, according to Bernstein inequality, $ \|
\nabla \Delta_j  W\|_2  \ge 2^{j-1}\|
 \Delta_j  W\|_2 $  and divide \re{n2} by $ \|\Delta_j W\|_2 $.
Then, estimating  the  commutator term thanks to \re{commu} and the last
term  thanks to \re{product} with $ s_1=s-1 $ and $ s_2=N/2 $,
using  Bernstein inequalities, it follows that
\begin{eqnarray}
\Re \frac{d}{dt} \|\Delta_j W\|_2  +
 \frac{(1-\omega)}{2} 2^{2j}
\|\Delta_j W\|_2 & \le & 2\|\Delta \Delta_j v\|_{L^2}+2  \|\Delta_j \div \tau\|_2  \nonumber \\
 &  &\hspace*{-34mm} +\gamma_j
 2^{-j(s-1)} (\| u\|_{B^{N/2+1}}+\| v\|_{B^{N/2+1}}) \|W\|_{B^{s-1}}
 \, , \label{n3}
\end{eqnarray}
with $ \|(\gamma_j) \|_{l^1(\Z)} \lesssim 1 $. Multiplying by
$2^{j(s-1)}$ and  summing in $ j  \in K_\varepsilon
 $, it follows that
\begin{eqnarray}
\Re \frac{d}{dt} \| W\|_{B^{s-1}_{K_\varepsilon}} & + &
\frac{(1-\omega)}{2} \| W\|_{ B^{s+1}_{K_\varepsilon}} -2  \|
\tau\|_{ B^{s}_{K_\varepsilon}}\nonumber \\
 & \le&  \quad
 (\| u\|_{ B^{N/2+1}}+\| v\|_{ B^{N/2+1}}) \|W\|_{  B^{s-1}}
 \, . \label{n4}
\end{eqnarray}
Proceeding in the same way with \re{u2} but summing in $ j \in
I\cup J_\varepsilon$,
 we obtain
 \begin{eqnarray}
\Re \frac{d}{dt} \| W\|_{B^{s-1}_{I\cup J_\varepsilon}} & + &
\frac{1}{2}\| W\|_{ B^{s+1}_{I\cup J_\varepsilon}} - \|Z\|_{
B^{s}_{I\cup
J_\varepsilon}}\nonumber \\
 & \le&  \quad
 (\| u\|_{B^{N/2+1}}+\| v\|_{ B^{N/2+1}}) \|W\|_{ B^{s-1}}
 \, . \label{n5}
\end{eqnarray}
Now, for $ j\in \Z $, we infer from \re{nonnewtonian} that
\begin{eqnarray}
& \varepsilon &  \partial_t \Delta_j \tau+\varepsilon
(u.\nabla)\Delta_j \tau  + \Delta_j \tau =2\omega \Delta_j
\DD[u]\nonumber \\
& & -  \varepsilon [(u.\nabla),\Delta_j]\tau +\varepsilon \Delta_j
\ga (\nabla u,\tau) \, . \label{m1}
\end{eqnarray}
Rewriting $ \Delta_j ( \tau-2\omega \DD[u]) $ as $ \Delta_j
(\tau-2\omega \DD[W]-2\omega \DD[v]) $
 for $ j\in  K_\varepsilon  $ and as
$$
\varepsilon^\beta \Delta_j \tau -2\omega\varepsilon^\beta \Delta_j (\DD[W]+\DD[v])+(1-\varepsilon^\beta) \Delta_j Z $$
 for $ j\in I\cup J_\varepsilon $,
 similar considerations as above lead to the two following inequalities
\begin{eqnarray}
\varepsilon\frac{d}{dt} \| \tau
\|_{B^{s}_{K_\varepsilon}} & + & \| \tau\|_{
B^{s}_{K_\varepsilon}} \le
  4\omega  \|W\|_{ B^{s+1}_{K_\varepsilon}}
 \nonumber \\
   & \ & +4\omega  \|v\|_{ B^{s+1}_{K_\varepsilon}}
   +C \, \varepsilon \, \|u\|_{B^{N/2+1}} \|\tau\|_{B^s}  \, ,\label{tau1Bs}
\end{eqnarray}
and
 \begin{eqnarray}
\varepsilon\frac{d}{dt} \| \tau \|_{B^{s}_{I\cup J_\varepsilon}}
+ \varepsilon^\beta \| \tau\|_{B^{s}_{I\cup J_\varepsilon}} & \le
&  \| Z \|_{ B^{s}_{I\cup J_\varepsilon}}
  + 4\omega \varepsilon^\beta \|W\|_{B^{s+1}_{I\cup
J_\varepsilon } }\nonumber \\
& & \hspace*{-10mm}+4\omega \varepsilon^\beta
\|v\|_{B^{s+1}_{I\cup J_\varepsilon}}
 + C \, \varepsilon \, \|u\|_{B^{N/2+1}}
\|\tau\|_{B^s} \, . \label{tau2}
\end{eqnarray}
 Adding $\re{n4}+\re{n5}+4 (\re{tau2}+   \re{tau1Bs}) $, \re{utaubesov} follows.
\subsection{Estimate on $ \tau -2\omega \DD[u]$}
\begin{lemma}
\begin{eqnarray}
\frac{d}{dt} \|Z\|_{B^{s-2}_{J_\varepsilon}} +
\frac{1}{\varepsilon} \|Z\|_{B^{s-2}_{J_\varepsilon}} & \le &
\frac{(1+\omega)}{\Re}  \|\tau\|_{B^s_{J_\varepsilon}}
+ \|{\cal P}f \|_{B^{s-1}_{J_\varepsilon}} \nonumber \\
 & & \hspace*{-25mm}+C \,\alpha \ln (\varepsilon^{-1})
 \Bigl(\|u\|_{B^{N/2+1}} +\|\tau\|_{B^{N/2}}\Bigr) \|u\|_{B^{s-1}}
 \label{ZZJ}
\end{eqnarray}
\begin{eqnarray}
\frac{d}{dt} \|Z\|_{B^{s}_{I}} + \frac{1}{\varepsilon}
\|Z\|_{B^{s}_{I}} & \le & \frac{(1+\omega)}{\Re}
\|\tau\|_{B^s_{I}}
+ \|{\cal P}f \|_{B^{s-1}_{I}} \nonumber \\
 & & \hspace*{-25mm}+C
 \Bigl(\|u\|_{B^{N/2+1}} +\|\tau\|_{B^{N/2}}\Bigr) \|u\|_{B^{s-1}}\, .
 \label{ZZI}
\end{eqnarray}
\end{lemma}

{\it Proof.}
Applying $ \Delta_j $ to \re{Z1} and taking the $ L^2 $-scalar
product with $ \Delta_j Z $ we get
$$
\frac{d}{dt}\|\Delta_j Z \|_{L^2}+\frac{(1-\omega)}{2\Re} 2^{2 j
}\|\Delta_j Z \|_{L^2}+ \frac{1}{\varepsilon} \|\Delta_j Z
\|_{L^2} \lesssim \|\Delta_j f_1\|_{L^2}+\|\Delta_j f_2\|_{L^2}\,
,
$$
 where
 $$
 f_1=\frac{2\omega}{\Re} \DD[{\cal P} \div \tau] -\frac{(1-\omega)}{\Re} \Delta \tau
  +\frac{2\omega}{\Re} \DD[{\cal P} f]
 -2\omega \DD[{\cal P}(u.\nabla)u]
 $$
 and
 $$ f_2= {\cal P}(u.\nabla)\tau +g(\nabla u,\tau) \quad .
 $$
Multiplying this inequality by $ 2^{j(s-2)} $, summing in $ j \in
J_\varepsilon
 $  we infer that
\begin{eqnarray}
\frac{d}{dt}\| Z\|_{B^{s-2}_{J_\varepsilon}} & + &
\frac{(1-\omega)}{2\Re} \| Z \|_{B^s_{J_\varepsilon}}+
\frac{1}{\varepsilon} \| Z \|_{B^{s-2}_{J_\varepsilon}} \le
\frac{(1+\omega)}{\Re}
 \|\tau \|_{B^s_{J_\varepsilon}} + \|{\cal P} f \|_{B^{s-1}_{J_\varepsilon}} \nonumber \\
&  & + \|(u.\nabla) u \|_{B^{s-1}_{J_\varepsilon}}+ \|(u.\nabla)
\tau \|_{B^{s-2}_{J_\varepsilon}}+\|\ga (\nabla u,\tau)
\|_{B^{s-2}_{J_\varepsilon}} \, . \label{n11}
\end{eqnarray}
For $ s>1 $ we estimate the nonlinear term thanks to \re{product}
with respectively $ (s_1,s_2)=(s-1,N/2), \, (s-1,N/2-1) $ and $
(s-2, N/2) $. For $ s=1 $ ( of course $ N=2$) we estimate the
first nonlinear term in the same way and use the following lemma
 to estimate the two last ones. This lemma follows directly from the
 definitions of $ I $ and $ J_\varepsilon $ and the fact  that, for $ |s|\le N/2 $,
 the usual product maps continuously\footnote{For $ 1\le p\le \infty$, $\| f \|_{B^{s,p}}
 =\left\| \{ 2^{js}\| \Delta_j (f) \|_{L^2} \} \right\|_{l^p(\Z)} \; .
$}
  $ B^{-s,1}\times B^{s,1} $ into $ B^{-N/2,\infty} $ (see for instance \cite{Taylor}). Note, in particular,
 that $ |J_\varepsilon|\lesssim \alpha \ln (\varepsilon^{-1})$.
\begin{lemma} \label{lm}
For all $ s_1, s_2 \le N/2$ with $ s_1+s_2=0 $ it holds
\begin{equation}
\|a\, b \|_{B^{-N/2}_{J_\varepsilon}}\lesssim \alpha \ln (
\varepsilon^{-1}) \|a\|_{B^{s_1}} \|b\|_{B^{s_2}}
\label{ko1}\quad .
\end{equation}
and
\begin{equation}
\|a\, b \|_{B^{-N/2+2}_{I}}\lesssim  \|a\|_{B^{s_1}}
\|b\|_{B^{s_2}} \label{ko2}\, .
\end{equation}
\end{lemma}

We apply this lemma with $ (s_1,s_2)=(0,0) $ and $ (-1,1) $ for respectively the second
and the third nonlinear term of \re{n11} to complete
 the proof of \re{ZZJ}. Finally \re{ZZI} can be easily obtained in
 the same way by using that $ \|a\|_{B^s_{I}} \le
 \|a\|_{B^{s'}_{I}} $ for $ s'\le s $ and \re{ko2}.
\subsection{Convergence to the Newtonian flow}
From now on we set $ \gamma(\omega)=(1-\omega)/2-16\omega
$ and  assume that $ 0\le \omega\le \omega_0 $ with $\gamma(\omega_0)>0 $. \\
We proceed as in Section 2.3. For $ 0<\beta<1 $, we add
\re{utaubesov} and $\varepsilon^{2\beta} (\re{ZZJ}+\re{ZZI})$ to get
\begin{eqnarray}
&  \displaystyle\frac{d}{dt}  & \Bigl( \Re \| W\|_{B^{s-1}}
+ 4\varepsilon  \|\tau\|_{B^{s}}
+\varepsilon^{2\beta}(\|Z\|_{B^{s-2}_{J_\varepsilon}}
+\|Z\|_{B^s_{I}}) \Bigr) \nonumber \\
& & + \gamma(\omega_0) \| W\|_{ B^{s+1}_{K_\varepsilon}} +\| W\|_{
B^{s+1}_{I\cup J_\varepsilon}}
 +2  \| \tau\|_{ B^{s}_{K_\varepsilon}}+
  \varepsilon^\beta
 \|\tau\|_{B^s_{I\cup
J_\varepsilon}}+
\frac{\varepsilon^{2\beta-1}}{2}\Bigl(\|Z\|_{B^{s-2}_{J_\varepsilon}}
+\|Z\|_{B^s_{I}} \Bigr)\nonumber \\
& \le & 16\omega \|v\|_{ B^{s+1}_{K_\varepsilon}}+16\omega
\varepsilon^\beta
  \|v\|_{B^{s+1}_{I\cup
J_\varepsilon}}+\varepsilon^{2\beta}\|{\cal P} f
\|_{B^{s-1}_{I\cup J_\varepsilon}}
   \nonumber \\
 & & + C \, \varepsilon \|u\|_{B^{N/2+1}} \|\tau\|_{B^s}+ C \,
 (\| u\|_{ B^{N/2+1}}+\| v\|_{ B^{N/2+1}}) \|W\|_{
 B^{s-1}}\nonumber \\
 & &  + C \,\alpha  \varepsilon^{2\beta} \ln(\varepsilon^{-1})
 ( \|u\|_{B^{N/2+1}} +\|\tau\|_{B^{N/2}})\, \|v\|_{B^{s-1}}\, .
 \label{utaubesov2}
\end{eqnarray}
Here we used that for $ \varepsilon $ small enough, $
\varepsilon^\beta\le \min(16\gamma(\omega_0), \frac{\Re}{4})
$, $ \varepsilon^{2\beta-1}/2\ge 5  $ and
$$
5\|Z\|_{B^s_{J_\varepsilon}} \lesssim
\varepsilon^{-2\alpha}\|Z\|_{B^{s-2}_{J_\varepsilon}}\le
\frac{\varepsilon^{2\beta-1}}{2} \|Z\|_{B^{s-2}_{J_\varepsilon}}
$$
as soon as
\begin{equation}
\label{condbetabesov} 0<\alpha<1/2 \quad \mbox{ and } \quad
0<2\beta<1-2\alpha \, .
\end{equation}
From now on we set $(\alpha,\beta)=(1/8,1/8) $ so that
\re{condbetabesov} is satisfied. Setting
\begin{eqnarray*}
X_s(t) & =& \Re \|W(t)\|_{B^{s-1}}+ 4\varepsilon \|\tau
\|_{B^s} +\varepsilon^{2\beta}
\|Z\|_{B^{s-2}_{I\cup J_\varepsilon}} \nonumber \\
 &  & \hspace*{-1cm}+\int_0^t \frac{\gamma(\omega_0)}{2} \|W\|_{B^{s+1}}+
 \Bigl(\|\tau\|_{B^s_{K_\varepsilon}}+\varepsilon^\beta
 \|\tau\|_{B^s_{I\cup J_\varepsilon }} \Bigr)
 + \frac{\varepsilon^{2\beta-1}}{2}\Bigl(
 \|Z\|_{B^{s-2}_{J_\varepsilon}}+
 \|Z\|_{B^{s}_{I}}\Bigr)\, ds \, ,
 \end{eqnarray*}
 we infer that
 \begin{eqnarray*}
 \frac{d}{dt} X_s(t) & \le &
 16\omega  \|v\|_{ B^{s+1}_{K_\varepsilon}}+16\omega
\varepsilon^\beta
  \|v\|_{B^{s+1}_{I\cup
J_\varepsilon}}+\varepsilon^{2\beta}\|{\cal P} f
\|_{B^{s-1}_{I\cup J_\varepsilon}}
   \nonumber \\
 & & + C(\Re,\omega)\,\Bigl[   \|W\|_{B^{N/2+1}} +\| v\|_{ B^{N/2+1}}
 +\varepsilon^{2\beta} \ln(\varepsilon^{-1})
 \|\tau\|_{B^{N/2}}\Bigr] X_s
 \nonumber \\
 & &  + C \,\alpha  \varepsilon^{2\beta} \ln(\varepsilon^{-1})
 ( \|v\|_{B^{N/2+1}} + \|W\|_{B^{N/2+1}}+\|\tau\|_{B^{N/2}})\, \|v\|_{B^{s-1}}\, .
 \end{eqnarray*}
 By Gronwall lemma we infer that
 \begin{eqnarray}
 X_s(t) & \le & \exp\Bigl( C(\omega,\Re) \Bigl(  \|v\|_{L^1_t B^{N/2+1}}
 + X_{N/2}(t) \Bigr)\Bigr)\nonumber \\
 &  & \hspace*{-5mm}\Bigl[X_s(0)+  16\omega  \|v\|_{L^1_t
B^{s+1}_{K_\varepsilon}}+16\omega  \varepsilon^\beta
  \|v\|_{L^1_t B^{s+1}_{I\cup
J_\varepsilon}}+\varepsilon^{2\beta}\|{\cal P} f \|_{L^1_t
B^{s-1}_{I\cup J_\varepsilon}}
   \nonumber \\
   & &\hspace*{-5mm} +C\, \alpha \ln(\varepsilon^{-1})\varepsilon^\beta \Bigl(
   X_{N/2}(t)\|v\|_{L^\infty_t B^{s-1}}+\varepsilon^\beta
\|v\|_{L^1_t B^{s-1}}\|v\|_{L^\infty_t B^{s-1}}\Bigr) \Bigr] \, .
\label{fo1}
 \end{eqnarray}
 where
 \begin{equation}
 X_s(0)=4\varepsilon
 \|\tau_0\|_{B^s}+\varepsilon^{2\beta}\Bigl(\|\tau_0-2\omega
 \DD[u_0]\|_{B^s_{I}}+\|\tau_0-2\omega
 \DD[u_0]\|_{B^{s-2}_{J_\varepsilon}}\Bigr) \, . \label{fo2}
 \end{equation}
Assuming that  $ T_\varepsilon^*\le T_0 $ and noticing that
 \begin{equation}
 \|v\|_{L^1_{T_0} B^{N/2+1}_{K_\varepsilon}} \to 0
 \mbox{ as } \varepsilon \to 0 \, , \label{fo3}
 \end{equation}
 we deduce from \re{fo1}-\re{fo2} and the continuity of $ t\mapsto
 X_{N/2}(t) $ that there exists $
 \varepsilon_0= \varepsilon_0(N,\|\tau_0\|_{B^{N/2}},{\cal P} f ,
u_0) $ such that for any
 $ 0<\varepsilon<\varepsilon_0 $ and any $ 0<t<T_\varepsilon^* $,
 $$
 X_{N/2}(t) \le \Lambda(\varepsilon)
 $$ with $ \Lambda(\varepsilon) \searrow 0 $ as $ \varepsilon \to 0
 $.
 This contradicts \re{exploseBs} of Theorem \ref{CMBesov} and thus  ensures that
 $ T_\varepsilon^*>T_0 $. \re{convuBs} and \re{epsilontauBs} follow as well.
To prove \re{epsilontauBs} we notice that from this last limit and
\re{fo3}, $ \|\tau-2\omega\DD[u]\|_{L^1_{T_0}
B^{N/2}_{K_\varepsilon}}\to 0 $,
$ \varepsilon^{2\beta-1} \|\tau-2\omega\DD[u]\|_{L^1_{T_0} B^{N/2}_{I}}\to 0 $ and $
\varepsilon^{2\beta-1}\|\tau-2\omega\DD[u]\|_{L^1_{T_0}
B^{N/2-2}_{J_\varepsilon}}\to 0 $. This gives the result since $
2\beta-1+2\alpha \le 0 $ and thus
$$
\|\tau-2\omega\DD[u]\|_{L^1_{T_0} B^s_{J_\varepsilon}}\lesssim
\varepsilon^{-2\alpha}\|\tau-2\omega\DD[u]\|_{L^1_{T_0}
B^{s-2}_{J_\varepsilon}}\lesssim
\varepsilon^{2\beta-1}\|\tau-2\omega\DD[u]\|_{L^1_{T_0}
B^{s-2}_{J_\varepsilon}}\, .
$$
Finally, for $ s>N/2 $, the proof follows the same lines using
that
 \begin{equation}
 \|v\|_{L^1_{T_0} B^{N/2+1}_{K_\varepsilon}} \le
 \varepsilon^{\alpha(s-N/2)}\|v\|_{L^1_{T_0}
 B^{s+1}_{K_\varepsilon}} \label{hoho}
 \end{equation}
 and thus with
 $$
 \varepsilon_0=\varepsilon_0(N,\|\tau_0\|_{B^{N/2}},
 \|u_0\|_{B^{N/2-1}}, \|v\|_{L^1_{T_0} B^{s+1}_{K_\varepsilon}}, \|{\cal P} f\|_{L^1_{T_0}
 B^{N/2-1}})  \, . $$
 This completes the proof of Theorems \ref{theoBs} and  \ref{critic}. \vspace*{2mm}\\
\section*{Acknowledgments}
 The authors are grateful to the Referees
for
 useful remarks.

\end{document}